\documentstyle[12pt]{amsart}

\newcommand{\kahler}{K\"ahler\ }
\newcommand{\wt}{\widetilde}
\newcommand{\PP}{{\mathbb P}}
\newcommand{\R}{{\mathbb R}}
\newcommand{\C}{{\mathbb C}}

\newcommand{\Z}{{\mathbb Z}}
\newcommand{\dbar}{\bar\partial}
\newcommand{\ddbar}{\partial\dbar}
\newcommand{\U}{{\rm U}}
\newcommand{\Oo}{{\cal O}}
\newcommand{\half}{{\frac{1}{2}}}
\newcommand{\vol}{{\operatorname{Vol}}}
\newcommand{\SU}{{\operatorname{SU}}}
\newcommand{\FS}{{{\operatorname{FS}}}}
\renewcommand{\phi}{\varphi}
\newcommand{\go}{\mathfrak}

\newtheorem{theo}{{\sc Theorem}}[section]
\newtheorem{cor}[theo]{{\sc Corollary}}
\newtheorem{lem}[theo]{{\sc Lemma}}
\newtheorem{prop}[theo]{{\sc Proposition}}
\newenvironment{proof}{\noindent{\em Proof:\/}}{\qed \medskip}

\title[Distribution
of zeros of sections of positive line bundles]{Distribution
of zeros of random and quantum chaotic sections of positive line bundles}

\author{Bernard Shiffman}
\author{Steve Zelditch}
\address{Department of Mathematics, Johns Hopkins University, Baltimore, MD
21218, USA} 
\email{shiffman@@math.jhu.edu {\it (first author),}
zel@@math.jhu.edu {\it (second author)}}

\thanks{Research of the first author
partially supported by  NSF grant \#DMS-9500491; research of the second
author partially supported by  NSF grant \#DMS-9703775.} 

\date{March 6, 1998}

\begin{document}

\begin{abstract}  We study the limit distribution of zeros of certain
sequences of holomorphic sections  of high powers $L^N$ of a positive
holomorphic Hermitian line bundle $L$ over a compact complex manifold $M$.
Our first result concerns `random' sequences of sections.  Using the
natural probability measure on the space of sequences of orthonormal bases
$\{S^N_j\}$ of $H^0(M, L^N)$, we show that for almost every sequence
$\{S^N_j\}$, the associated sequence of zero currents $\frac{1}{N}
Z_{S^N_j}$ tends to the curvature form $\omega$ of $L$.  Thus, the zeros of
a  sequence of sections $s_N \in H^0(M, L^N)$ chosen independently and at
random become uniformly distributed. Our second result concerns the zeros
of quantum ergodic eigenfunctions, where the relevant orthonormal bases
$\{S^N_j\}$ of $H^0(M, L^N)$ consist of eigensections of a quantum
ergodic map. We show that also in this case the zeros become
uniformly distributed. \end{abstract}

\maketitle

\section{Introduction} 

This paper is concerned with the limit distribution of zeros of `random'
holomorphic sections and of `quantum ergodic' eigensections of powers of a
positive holomorphic line bundle $L$ over a compact complex manifold $M$.
To introduce our subject, let us consider the simplest case where $M=\C\PP^m$
and where $L$ is the hyperplane section bundle.  As
is well-known, sections of $L^{N}$ are given by homogeneous
polynomials $p_N(z_0, z_1, \dots, z_m)$ of degree $N$ on $\C^{m+1}$;  these
polynomials are called $\SU(m+1)$ polynomials when we consider them as
elements of a measure space with an $\SU(m+1)$-invariant Gaussian measure (see
\S\ref{section-SU}).  We are concerned with the question: what is the
limit distribution of zeros $Z_N = \{p_N = 0\} \subset M$ of a sequence
$\{p_N\}$ of such polynomials as the degree $N \rightarrow \infty$?   Of
course, if we consider all possible sequences, then little can be said.
However, if we  consider only the typical behavior, then there is a simple
answer: if the sequence $\{p_N\}$ is chosen independently and at random from
the ensembles of  homogeneous polynomials of degree $N$ and ${\cal
L}^2$-norm one, then the zero sets of $\{p_N\}$ almost surely become
uniformly distributed with respect to the volume form induced
by $\omega.$  

The same conclusion is true for any positive Hermitian holomorphic line
bundle $(L,h)$ over any compact complex manifold $M$. In place of
homogeneous polynomials of degree $N$, one now considers holomorphic
sections $s_N \in H^0(M, L^{N})$.  The curvature form $\omega=c_1(h)$ of $h$
defines a \kahler structure on $M$, and the metrics $h,\;\omega$ provide a
Hermitian inner product on $H^0(M, L^{N})$. (See equations
(\ref{dV})--(\ref{inner}) in \S \ref{background}.) We then have the notion of
a `random' sequence of ${\cal L}^2$-normalized sections of $H^0(M, L^{N})$.
Namely, we consider the probability space $({\cal S}, d\mu)$, where ${\cal
S}$ equals the product $\prod_{N=1}^{\infty} SH^0(M, L^{N})$ of the unit
spheres $SH^0(M, L^{N})$ in $H^0(M, L^{N})$ and $\mu$ is the product of Haar
measures on these spheres. Given a sequence ${\bf s} = \{s_N\}\in {\cal S}$,
we associate the currents of integration $Z_{s_N}$ over the zero divisors of
the sections $s_N$. In complex dimension 1, $Z_{s_N}$ is simply the sum of
delta functions at the zeros of $s_N$.  Our first result states that for a
random (i.e., for almost all) ${\bf s}\in{\cal S}$, the sequence of zeros of
the sections $s_N$ are asymptotically uniformly distributed:

\begin{theo} \label{zerotheorem} For $\mu$-almost all ${\bf
s}=\{s_N\}\in{\cal S}$, $\frac{1}{N}Z_{s_N}\to \omega$ weakly in the sense
of measures; in other words, $$\lim_{N\to\infty}\left(\frac{1}{N}
Z_{s_N},\phi\right) = \int_M\omega\wedge\phi$$ for all continuous
$(m-1,m-1)$ forms $\phi$. In particular, $$\lim_{N \rightarrow
\infty}\frac{1}{N} \vol_{2m-2}\{z\in U:s_N(z)=0\}=m\vol_{2n}U\;,$$ for
$U$ open in $M$ (where $\vol_k$ denotes the Riemannian $k$-volume in
$(M,\omega)\,$).\end{theo}

The key ideas in the proof of Theorem~\ref{zerotheorem} (as well as
Theorem~\ref{RANTH} below) are Tian's theorem \cite{Ti,Z4} on approximating
the metric $\omega$ using the sections of $H^0(M, L^{N})$ 
(see Theorem~\ref{tianstheorem}) and an asymptotic estimate of the variances
of $Z_{s_N}$, regarded as a current-valued random variable
(Lemma~\ref{variance}).

A closely related issue is the distribution of zeros of sections
$\{S_j^N\}$ forming random orthonormal bases of $H^0(M, L^{N})$.
Such bases are increasingly used to model orthonormal bases of  quantum
chaotic eigenfunctions; e.g., see \cite{B.B.L, H,LeB.S, N.V}.  The
properties of these bases are very similar to those of
random orthonormal bases of spherical harmonics studied in \cite{Z1} and
\cite{V}.  To study the zeros of random orthonormal bases, we introduce the
probability space $({\cal ONB},d\nu)$, where ${\cal ONB}$ is the infinite
product of the sets ${\cal ONB}_N$ of orthonormal bases of the spaces
$H^0(M, L^{N})$, and $\nu =\prod_{N=1}^{\infty} \nu_N$, where $\nu_N$ is
Haar probability measure on ${\cal ONB}_N$. A point of ${\cal ONB}$ is thus
a sequence ${\bf S} = \{(S^N_1, \dots, S^N_{d_N}) \}_{N\geq 1}$ of
orthonormal bases (where $d_N=\dim H^0(M, L^{N})$), and we may ask whether
all of the zero sets $Z_{S^N_j}$ are tending simultaneously to the uniform
distribution.  The answer is still essentially yes,  but for technical
reasons we have to delete a subsequence of relative density zero of the
sections. 

\begin{theo} \label{RANTH} For $\nu$-almost all ${\bf S} = \{(S^N_1, \dots,
S^N_{d_N}) \} \in {\cal ONB}$, we have $$\frac{1}{d_N} \sum_{j=1}^{d_N}
\left( \frac{1}{N} Z_{S_j^N} - \omega, \phi \right)^2 \rightarrow 0$$  for
all continuous $(m-1,m-1)$ forms $\phi$.
Equivalently, for each $N$ there exists a subset $\Lambda_N \subset \{ 1,
\dots, d_N\}$ such that $\frac{\# \Lambda_N}{d_N} \rightarrow 1$ and
$$\lim_{N \rightarrow \infty, j \in \Lambda_N}\frac{1}{N}  Z_{S^N_j} =
\omega$$ weakly in the sense of measures.  \end{theo}

Our final result pertains to actual quantum ergodic eigenfunctions rather
than to random sections and shows that their zero divisors also become
uniformly distributed in the high power limit.  Recall that a quantum map is
a unitary operator which `quantizes' a symplectic map on a symplectic
manifold.  In our setting, the symplectic manifold is the \kahler manifold
$(M,\omega)$ and the map is a symplectic transformation
$\chi : (M, \omega) \rightarrow (M, \omega)$.  Under certain conditions,
$\chi$ may be quantized  as a sequence of unitary operators $U_{\chi, N}$ on
$ H^0(M, L^{N})$. The sequence defines a semiclassical Fourier
integral operator of Hermite type (or equivalently a semiclassical Toeplitz
operator). For the precise definitions and conditions, we refer to 
\cite{Z3}.  We call $U_{\chi, N}$ a `quantum ergodic map' if $\chi$ is also
an ergodic transformation of $(M, \omega)$.

\begin{theo}\label{QE} Let $(L, h) \rightarrow (M, \omega)$ be a positive
Hermitian line bundle over a \kahler manifold with $c_1(h) = \omega$ and let
$U_{\chi, N} : H^0(M, L^{N}) \mapsto  H^0(M, L^{N})$ be a
quantum ergodic map.  Further, let  $\{S_1^N, \dots, S_{d_N}^N\}$ be an
orthonormal basis of eigensections of $U_{\chi, N}$.  
Then there exists a subsequence $\Lambda \subset \{(N,j): N = 1, 2, 3,
\dots, j \in  \{1, \dots, d_N\} \}$ of density one such that $$\lim_{N
\rightarrow \infty, (N,j) \in \Lambda} \frac{1}{N}Z_{S^N_j} =
\omega$$ weakly in the sense of measures. \end{theo} 
 
This result was proved  independently  by
Nonnenmacher-Voros \cite{N.V} in the case of the theta bundle over an
elliptic curve $\C / \Z^2$.  The main step is to establish the following
result:   

\begin{lem} Let $(L,h) \rightarrow (M, \omega)$ be a positive
Hermitian holomorphic line bundle over a \kahler manifold $M$ with 
$c_1(h) = \omega$. Let $s_N \in H^0(M, L^{N})$, $N=1,2,\ldots$, be
a sequence of  sections with the property that $\|s_N(z)\|^2 \rightarrow 1$
in the weak* sense as $N \rightarrow \infty$.  Then
$\frac{1}{N}Z_N\to\omega$ weakly in the sense of measures.
\label{W*}\label{mainlemma}\end{lem}

The convergence hypothesis means that $\int_M
\phi(z)\|s_N(z)\|^2dz\to\int_M\phi(z)dz$ for all $\phi\in{\cal C}^0(M)$. Our
proof of Lemma~\ref{mainlemma} is somewhat different and more general than
that of \cite{N.V}, but both are based on potential theory.  The
lemma was motivated by an analogous result of Sodin
\cite{So} on the asymptotic equidistribution of zero sets of sequences of
rational functions  in one variable (see also \cite{R.Sh,R.S} for the 
higher dimensional case); Sodin's result in turn arose from the
Brolin-Lyubich Theorem in complex dynamics (cf., \cite{F.S}). The
connection between Lemma~\ref{mainlemma} and Theorems \ref{RANTH}, \ref{QE}
will be established in \S\ref{section-ergodic}, the main point being
that both random orthonormal bases and orthonormal bases of chaotic
eigenfunctions satisfy the hypothesis of the lemma 
(Theorems~\ref{EONB}, \ref{Shn}).

We end this introduction with a brief discussion of related results.  There
is an extensive literature on the distribution of zeros of random
polynomials, beginning with the classical papers of Bloch-Polya \cite{B.P},
Littlewood-Offord \cite{L.O}, Kac \cite{K.1} and Erdos-Turan \cite{E.T} on
polynomials in one variable. The articles of Bleher-Di \cite{Bl.D} and
Shepp-Vanderbei \cite{S.V} contain recent results and further references.  In
addition to the mathematical literature there is a growing physics literature
on zeros of random polynomials and chaotic quantum eigenfunctions, see in
particular \cite{Bl.D,B.B.L,H,LeB.S,N.V}.  As in this paper, these articles
are largely concerned with the distribution of zeros in the semiclassical
limit.  The main theme is that the distribution of zeros of eigenfunctions
of quantum maps should reflect the signature of the dynamics of the
underlying classical system: in the case of ergodic quantum maps, the zeros
should be uniformly distributed in the semiclassical limit while in the
completely integrable case they should concentrate in a singular way. 
Random polynomials (or more generally sections) are believed to provide an
accurate model for quantum chaotic eigenfunctions and hence there is 
interest in understanding how their  zeros are distributed and how the zeros
are  correlated.

To our knowledge, the prior  results on distribution of zeros of random
holomorphic sections only go as far as determining the average
distribution.  In the special case of SU(2) polynomials it is shown in 
\cite{B.B.L} that the average distribution is uniform. Our result that the
expected distribution is achieved asymptotically by almost every sequence of
sections appears to be new even in that case.   Regarding zeros of quantum
ergodic eigenfunctions, the only prior rigorous result appears to be that of
\cite{N.V} mentioned above. We should also mention the study of the zeros of
certain sections of positive line bundles in the almost complex setting 
which has recently been made by Donaldson \cite{Don};  the relevant zero sets
were also shown to be  uniformly distributed in the high power limit.

\smallskip\noindent{\it Acknowledgments:\/}  We would like to thank
S.~Nonnenmacher and A.~Voros for sending us a copy of their paper \cite{N.V}
prior to publication and to acknowledge their priority on the overlapping
result.  We would also like to thank W.~Minicozzi for discussions
of Donaldson's paper at the outset of this work and for suggesting that we
study random sequences of sections.

\section{Background}\label{background} We begin by introducing some
terminology and basic properties of orthonormal bases of holomorphic
sections of powers of a positive line bundle.

\subsection{Notation}\label{notation}

Throughout this paper, we let $L$ denote an ample holomorphic line bundle
over an $m$-dimensional compact complex (projective) manifold $M$. We
denote the space of global holomorphic sections of $L$ by $H^0(M,L)$.  We
let ${\cal D}^{p,q} (M)$ denote the space of ${\cal C}^\infty$
$(p,q)$-forms on $M$, and we let ${\cal D}'{}^{p,q}(M) = {\cal
D}^{m-p,m-q}(M)'$ denote the space of $(p,q)$-currents on $M$; $(T,\phi) =
T(\phi)$ denotes the pairing of $T\in {\cal D}'{}^{p,q}(M)$ and $\phi \in
{\cal D}^{m-p,m-q}(M)$.  If $L$ has a smooth Hermitian metric $h$, its {\it
curvature form} $c_1(h)\in {\cal D}^{1,1}(M)$ is given locally by
$$c_1(h)=-\frac{\sqrt{-1}}{\pi}\ddbar \log \|e_L\|_h\;,$$ where $e_L$ is a
nonvanishing local holomorphic section of $L$, and
$\|e_L\|_h=h(e_L,e_L)^{1/2}$ denotes the $h$-norm of $e_L$.  The curvature
form $c_1(h)$ is a de Rham representative of the Chern class $c_1(L)\in
H^2(M,\R)$; see \cite{G.H,S.S}.  Since $L$ is ample, we can give $L$ a
metric $h$ with strictly positive curvature form, and we give $M$ the
\kahler metric $\omega=c_1(h)$.  Then $\int_M \omega^m = c_1(L)^m \in \Z^+$.
Finally, we give $M$ the volume form \begin{equation}\label{dV} dV=
\frac{1}{c_1(L)^m} \omega^m\;,\end{equation} so that $M$ has unit volume:
$\int_M dV=1$.

This paper is concerned with the spaces $H^0(M,L^N)$ of sections of
$L^N=L\otimes\cdots\otimes L$.  The metric $h$ induces Hermitian metrics
$h_N$ on $L^N$ given by $\|s^{\otimes N}\|_{h_N}=\|s\|_h^N$.  We give
$H^0(M,L^N)$ the inner product structure
\begin{equation}\label{inner}\langle s_1, s_2 \rangle = \int_M h_N(s_1,
s_2)dV \quad\quad (s_1, s_2 \in H^0(M,L^N)\,)\;,\end{equation} and we write
$|s|=\langle s,s \rangle^{1/2}$.   We let $d_N=\dim H^0(M,L^N)$.  It is well
known that for $N$ sufficiently large, $d_N$ is given by the Hilbert
polynomial of $L$, whose leading term is $\frac{c_1(L)^m}{m!}N^m$ (see, for
example \cite[Chapter~7]{S.S}).

For a holomorphic section $s\in H^0(M, L^{N})$, we
let $Z_s$ denote the current of integration over the zero divisor of $s$. 
In a local frame $e_L^N$ for $L^{N}$, we can write $s=\psi e_L^N$, where
$\psi$ is a holomorphic function.  We recall the Poincar\'e-Lelong formula
\begin{equation} \label{PL} Z_s =\frac{\sqrt{-1}}{\pi} \ddbar \log|\psi| =
\frac{\sqrt{-1}}{\pi} \ddbar \log \|s\|_{h_n} + N\omega \;.\end{equation}
We also consider the normalized zero divisor $$\wt Z_s^N=\frac{1}{N}Z_s\;,$$
so that the currents $\wt Z_s^N$ are de Rham representatives of $c_1(L)$,
and thus \begin{equation} \left(\wt Z_s^N,\omega^{m-1}\right) =
\frac{c_1(L)^m}{m!}\;. \label{constantmass}\end{equation} Equation
(\ref{constantmass}) says that the currents $\wt Z_s^N$ all have the same
mass.

For example, we consider the hyperplane section bundle, denoted $\Oo(1)$, over
$\C\PP^m$.  Sections $s\in H^0(\C\PP^m,\Oo(1))$ are linear functions on
$\C^{m+1}$; The zero divisors $Z_s$ are projective hyperplanes. The line
bundle $\Oo(1)$ carries a natural metric $h_\FS$ given by
$$\|s\|_{h_\FS}([w])=\frac{|(s,w)|}{|w|}\;, \quad\quad
w=(w_0,\dots,w_m)\in\C^{m+1}\;,$$ for
$s\in\C^{m+1*}\equiv H^0(\C\PP^m,\Oo(1))$, where $|w|^2=\sum_{j=0}^m
|w_j|^2$ and $[w]\in\C\PP^m$ is the complex line through $w$. The
curvature form of $h_\FS$ is given by \begin{equation} c_1(h_\FS)= 
\omega_\FS=\frac{\sqrt{-1}}{2 \pi} \ddbar \log |w|^2 \;,\end{equation}
where $\omega_\FS$ is the Fubini-Study \kahler form on
$\C\PP^m$. Here, $\omega_\FS$ is normalized so that it represents the
generator of $H^2(\C\PP^m,\Z)$.  The $N$-th tensor power of $\Oo(1)$ is
denoted $\Oo(N)$.  Elements of $H^0(\C\PP^m,\Oo(N))$ are homogeneous
polynomials on $\C^{m+1}$ of degree $N$; hence, $\dim H^0(\C\PP^m,\Oo(N)) =
{N+m \choose m} = \frac{1}{m!}N^m + \cdots$.

\subsection{Holomorphic sections and CR holomorphic functions}

The setting for our analysis is the Hardy space $H^2(X) \subset {\cal
L}^2(X)$ where $X \rightarrow M$ is the principal $S^1$ bundle associated to
$L$.  To be precise, let $L^*$ be the dual line bundle to $L$ and let $D =
\{v \in L^* : h(v,v)< 1\}$ be its unit disc bundle relative to the metric
induced by $h$ and let $X=\partial D = \{v \in L^* : h(v,v)= 1\}$. The
positivity of $c_1(h)$ is equivalent to the disc bundle $D$ being strictly
pseudoconvex in $L^*$ (see \cite{Gr}).
 
We let $r_{\theta}x =e^{i\theta} x$ ($x\in X$) denote the $S^1$ action on
$X$  and denote its infinitesimal generator  by $\frac{\partial}{\partial
\theta}$.  As the boundary of a strictly pseudoconvex domain, $X$ is a CR
manifold, and the Hardy space $H^2(X)$ mentioned above is by definition the
space of square integrable CR functions on $X$.  Equivalently, it is the
space of boundary values of holomorphic functions on $D$ which are in
${\cal L}^2(X)$. The $S^1$ action on $X$ commutes with the
Cauchy-Riemann operator $\bar{\partial}_b$; hence $H^2(X) = \bigoplus_{N =
0}^{\infty} H^2_N(X)$ where $H^2_N(X) = \{ f \in H^2(X): f(r_{\theta}x) =
e^{i N \theta} f(x) \}.$  A section $s$ of $L$ determines an equivariant
function $\hat{s}$ on $L^*$ by the rule: $\hat{s}(z, \lambda) = 
\left(\lambda, s(z) \right)$ ($z \in M, \lambda \in L^*_z$). It is clear
that if $\tau \in \C$ then $\hat{s}(z, \tau \lambda) =  \tau \hat{s}.$  We
will usually restrict $\hat{s}$ to $X$ and then the equivariance property
takes the form: $\hat{s}(r_{\theta} x) = e^{i  \theta}\hat{s}(x).$ 
Similarly, a section $s_N$ of $L^{N}$ determines an equivariant function
$\hat{s}_N$ on $L^*$: put $\hat{s}_N(z, \lambda) = \left(
\lambda^{N}, s_N(z) \right)$ where $\lambda^{N} = \lambda
\otimes \cdots\otimes \lambda$; then $\hat s_N(r_\theta x) = e^{iN\theta}
s_N(x)$. The map $s \mapsto \hat{s}$ is a unitary equivalence between $H^0(M,
L^{N})$ and $H^2_N(X).$

We now recall the strong form of Tian's theorem \cite{Ti} given in
\cite{Z4}:

\begin{theo} {\rm \cite{Z4}} Let $M$ be a compact complex  
manifold of dimension $m$ (over $\C$) and let $(L, h) \rightarrow M$ be a
positive Hermitian holomorphic line bundle.  Let $\{S_1^N, \dots, S_{d_N}^N
\}$ be any orthonormal basis of $H^0(M, L^{N})$ (with respect to the inner
product defined above). Then there exists a complete asymptotic expansion $$
\sum_{j=1}^{d_N} \|S^N_j(z)\|_{h_N}^2 = a_0 N^m + a_1(z) N^{m-1} + a_2(z)
N^{m-2} + \dots$$ with $a_0 = \frac{c_1(L)^m}{m!}$ and with the lower
coefficients $a_j(z)$ given by invariant polynomials in the higher derivatives
of
$h$. More precisely, for any $k\geq 0$,
$$\Big\| \sum_{i=0}^{d_N} \|S^N_i\|_{h_N}^2 - \sum_{j < R} a_j
N^{m-j}\Big\|_{{\cal C}^k} \leq C_{R,k} N^{m- R}. \label{tianstheorem}
$$\end{theo}
 
Note that since the $S^N_j$ have unit length (as elements of $H^0(M, L^{N})$),
if we integrate the above asymptotic expansion over $M$ (with respect to the
volume $dV$, we get simply $d_N$.  Thus the integrals of the $a_j$ are the
coefficients of the Hilbert polynomial of $L$. (The constant $a_0$ differs
from that of \cite{Ti} and \cite{Z4}, since we use here the normalized
volume $dV$ on $M$.)

The canonical map
\begin{equation} {\Phi}_N : M  \rightarrow \PP H^0(M, L^{\otimes
N})^*,\;\;\;\;\; z \mapsto \{s \in H^0(M, L^{\otimes N}) : s(z) =
0\}\end{equation} can be described in terms of an orthonormal basis ${\bf
S}=\{S_1^N, \dots, S_{d_N}^N \}$ by the map 
\begin{equation} \label{kodaira} 
\Phi_N^{\bf S}: M \rightarrow \C\PP^{d_N-1},\;\;\;\;\; z \mapsto [S_1^N(z),
\dots, S_{d_N}^N(z)]\;.\end{equation} We shall drop the {\bf S} and denote
the map given in (\ref{kodaira}) simply by $\Phi_N$.  For $N$ sufficiently
large, the sections $\{S_1^N, \dots, S_{d_N}^N
\}$ do not have common zeros and (\ref{kodaira}) gives a holomorphic
embedding, by the Kodaira embedding theorem; see \cite{G.H, S.S}. 

Theorem~\ref{tianstheorem} can be regarded as an asymptotic formula for the
distortion function between the metrics $h_N$ and $\Phi_N^*h_\FS$ on the line
bundle $L^{N}$.  It also gives the following asymptotic estimate of the
Riemannian distortion of the maps $\Phi_N$:

\begin{cor} {\rm \cite{Z4}} Let $\omega_\FS$ denote the Fubini-Study form on
$\C  \PP^{d_N-1}$. Then for any $k\geq 0$, $$\left\|\frac{1}{N} 
\Phi_N^*(\omega_\FS) -  \omega\right\|_{{\cal C}^k} = O(\frac{1}{N})\;.$$
\label{tianscor} \end{cor}

\section{Zeros of random sections}

Our first aim is to determine the expected value of the normalized zero
divisor $\wt Z_s$ as $s$ is chosen at random from the unit sphere $$SH^0(M,
L^{N}):=\{s\in H^0(M, L^{N}):|s|=1\}$$ (or equivalently as $[s]\in\PP
H^0(M, L^{N})$ is chosen at random with respect to the Fubini-Study
volume). As above, we fix one orthonormal basis $\{S^N_j\}$ of $H^0(M,
L^{N})$ and write $S^N_j = f_j e_L^N$ relative to a holomorphic frame (=
nonvanishing section) $e_L^N$ over an open set $U\subset M$.  Any section
in $SH^0(M, L^{N})$ may then be written as $s = \sum_{j=1}^{d_N} a_j f_j
e_L^N$ with $\sum_{j=1}^{d_N} |a_j|^2 = 1.$ To simplify the notation we let
$f=(f_1,\ldots,f_{d_N}):U\to \C^{d_N}$ (which is a local representation of
$\Phi_N$) and we put $$ \sum_{j=1}^{d_N} a_j f_j = \langle a, f \rangle .$$
Hence \begin{equation}\label{PL2} \wt Z^N_s = \frac{\sqrt{-1}}{N\pi} \ddbar
\log|\langle a, f \rangle|\;.\end{equation}

\subsection{Expected distribution of zeros} \label{ezd}

We shall frequently use the notation $E(Y)$ for the expected value of a
random variable $Y$ on a probability space $(\Omega, d \mu)$, i.e. $E(Y) =
\int_{\Omega}Y d\mu$.

We view $\wt Z_s^N$ as a ${\cal D}'{}^{1,1}(M)$-valued random variable (which
we call simply a `random current') as $s$ varies over $SH^0(M, L^{N})$
regarded as a probability space with the standard measure, which we denote by
$\mu_N$.  The expected distribution of zeros of the random section $s$ is the
current $E (\wt Z_s^N)\in {\cal D}'{}^{1,1}(M)$ given by \begin{equation}
\left(E (\wt Z_s^N),\phi\right) = \int_{S^{2d_N-1}} \left(\wt
Z_s^N,\phi\right)d\mu_N , \quad \phi\in {\cal D}^{m-1,m-1}(M),\end{equation}
where we identify $SH^0(M, L^{N})$ with the unit $(2d_N+1)$-sphere
$S^{2d_N-1}\subset\C^{d_N}$. In fact, we have the following simple formula
for the expected zero-distribution in terms of the map $\Phi_N$ given by
equation (\ref{kodaira}):

\begin{lem} For $N$ sufficiently large so
that $\Phi_N$ is defined, we have:
$$E (\wt Z_s^N) =\frac{1}{N} \Phi_N^*\omega_\FS$$ \label{mean} 
\end{lem}

Lemma \ref{mean} is a special case of Lemma~\ref{mean2} below.  We give here a
short alternate proof of Lemma~\ref{mean} which will serve as an introduction
to our estimate on the variance (Lemma~\ref{variance}) to be given below. We
write \begin{equation} \omega_N = \frac{1}{N} \Phi_N^*\omega_\FS\;.
\end{equation}
In terms of our fixed
orthonormal basis, we have: \begin{equation} \omega_N =
\frac{\sqrt{-1}}{2 \pi N} \ddbar \log \sum_{j=1}^{d_N} |f^N_j|^2 =
\frac{\sqrt{-1}}{2 \pi N} \ddbar \log |f|^2,
\end{equation}
where $f=(f_0,\ldots,f_{d_N})$ is a local representation of $\Phi_N$ as 
defined above.
Let $\phi$ be a smooth $(m-1,
m-1)$ form, which we shall refer to as a `test form'. We may assume that we
have a coordinate frame for $L$ on Support$\;\phi$.  By (\ref{PL2}), we must
show that
\begin{equation}\frac{\sqrt{-1}}{\pi N} \int_{S^{2d_N-1}} \int_M \ddbar \log
|\langle a, f\rangle|\wedge \phi d\mu_N(a) = 
(\omega_N, \phi)\;.\end{equation}  To
compute the integral, we write $f = |f| u$ where $|u| \equiv 1.$ Evidently,
$\log |\langle a, f\rangle| = \log |f| + \log |\langle a, u \rangle|$. The
first term gives  \begin{equation} \frac{\sqrt{-1}}{\pi N} \int_M \ddbar \log
|f|\wedge \phi = \int_M \omega_N \wedge \phi. \end{equation} We now
look at the second term.  We have \begin{equation} \label{term=0} 
\frac{\sqrt{-1}}{\pi} \int_{S^{2d_N-1}} \int_M \ddbar \log |\langle a,
u\rangle|\wedge \phi d\mu_N(a)= \frac{\sqrt{-1}}{ \pi}\int_M  \ddbar \left[
\int_{S^{2d_N-1}} \log |\langle a, u \rangle|  d\mu_N(a)\right]\wedge \phi
=0,\end{equation} since the average $\int \log |\langle a, \omega \rangle| 
d\mu_N(a)$ is a constant independent of $u$ for $|u|=1$, and thus the
operator $ \ddbar$ kills it. \qed

Combining Corollary~\ref{tianscor} and Lemma~\ref{mean}, we obtain:

\begin{prop}  
$E(\wt Z^N_s) =\omega + O(\frac{1}{N})$; i.e., for each smooth test form
$\phi$, we have $$E(\wt Z_s^N,\phi)=\int_M \omega \wedge \phi
+O(\frac{1}{N})\;.$$ \label{mean1}\end{prop}

\subsection{Variance estimate}

The purpose of this section is to obtain the variance estimate we need to
obtain Theorem~\ref{RANTH}.  Let $\phi$ be a test form.  It follows from our
formula for the expectation (Lemma~\ref{mean}) that the variance of $(\wt
Z^N_s, \phi)$ is given by
\begin{equation}\label{def-var} E\left((\wt Z^N_s-\omega_N,\phi)^2\right)
=E\left(|(\wt Z^N_s,\phi)-(\omega_N,\phi)|^2\right)
=E\left((\wt Z^N_s,\phi)^2\right) - (\omega_N,\phi)^2\;. \end{equation}
We have the following estimate of the variance:

\begin{lem} \label{variance} Let
$\phi$ be any smooth test form.  Then $$E\left(|(\wt
Z^N_s,\phi)-(\omega_N,\phi)|^2\right) = O(\frac{1}{N^2}).$$ \end{lem}

\begin{proof} We again let $f$ be a local representation of
$\Phi_N$. Using (\ref{PL2}) we easily obtain

\begin{equation} E\left((\wt Z^N_s,\phi)^2\right) = \frac{-1}{\pi^2 N^2} 
\int_M \int_M (\ddbar \phi(z)) (\ddbar \phi(w)) \int_{S^{2d_N-1}} \log
|\langle f(z), a\rangle| \log |\langle f(w), a \rangle| d\mu_N(a) 
\end{equation} As in the previous lemma we write $f = |f| u$ with $|u|\equiv
1.$ Then \begin{eqnarray*}\log |\langle f(z), a\rangle| \log |\langle f(w),
a\rangle| &=& \log |f(z)| \log |f(w)| + \log|f(z)| \log |\langle u(w), a
\rangle| \\&&+ \log |f(w)| \log |\langle u(z), a \rangle | + 
 \log |\langle u(w), a \rangle |  \log |\langle u(z), a \rangle
|.\end{eqnarray*} 
The first term contributes 
\begin{equation}\frac{-1}{\pi^2 N^2}  \int_M \int_M (\ddbar \phi(z))
(\ddbar \phi(w))\log |f(z)| \log |f(w)| = \frac{1}{N^2}(\phi, \Phi_N^*
\omega_\FS)^2 =  (\phi, \omega_N )^2. \end{equation}
The middle two terms contribute zero to the integral by (\ref{term=0}). The
lemma at hand thus comes down to the following claim: 
\begin{equation} \label{claim}
\left| \int_M \int_M (\ddbar \phi(z)) (\ddbar \phi(w))
\int_{S^{2d_N-1}} \log |\langle u(z), a\rangle| \log |\langle u(w), a
\rangle| d\mu_N(a)\right|  = O(1).\end{equation} It suffices to show that
\begin{equation}\label{GN} G_N(x,y) := \int_{S^{2d_N-1}} \log |\langle x,
a\rangle| \log |\langle y, a \rangle| d\mu_N(a) = C_N+ O(1)\quad (x,y\in
S^{2d_N-1}), \end{equation} where $C_N$ is a constant and the $O(1)$ term is
uniformly bounded on $S^{2d_N-1} \times S^{2d_N-1}$. To verify (\ref{GN}), we
consider the Gaussian integral \begin{equation}\label{gaussian} \wt
G_N(x,y):=\int_{\C^{d_N}} e^{-|a|^2} \log |\langle x, a\rangle| \log
|\langle y, a \rangle| d a. \end{equation} We evaluate (\ref{gaussian}) in two
different ways. First, we use spherical coordinates $a = \rho \sigma$ with
$\sigma \in S^{2d_N-1}$. We have \begin{equation} \wt G_N(x,y) =
\int_0^{\infty} \int_{S^{2d_N-1}} e^{- \rho^2} \rho^{2d_N-1} \left(\log \rho +
\log|\langle x, \sigma \rangle|\right)\left(\log \rho + \log|\langle y, \sigma
\rangle|\right) d\rho d\sigma \end{equation} where $d\sigma$ denotes the
(non-normalized) volume element on the unit sphere.  Multiplying out we get
four terms. The only term that is non-constant is the term containing both $x$
and $y$. We then have \begin{eqnarray*}\wt G_N(x,y)&=& C_N +
\left[\int_0^{\infty} e^{- \rho^2} \rho^{2d_N-1} d\rho\right]
\int_{S^{2d_N-1}} \log|\langle x, \sigma \rangle| \log|\langle y, \sigma
\rangle| d\sigma \\& = & C_N+ \frac{(d_N-1)!}{2} \int_{S^{2d_N-1}}
\log|\langle x, \sigma \rangle| \log|\langle y, \sigma \rangle| d\sigma .
\end{eqnarray*}  We now evaluate $\wt G_N(x,y)$ 
a second way by noting that coordinates in $\C^{d_N}$ may be chosen so that
$x=(1,0,\ldots,0),\ 
y=(\zeta_1,\zeta_2,0,\ldots,0)$. Write $a'=(a_1,a_2),\
\tilde a=(a_3,\ldots,a_{d_N}),\ \zeta'=(\zeta_1,\zeta_2)$. Then
the integral becomes \begin{equation}\wt G_N(x,y) = \left[\int_{\C^{d_N-2} }
e^{-|\tilde{a}|^2} d\tilde{a}\right] \psi(\zeta') = \pi^{d_N-2} \psi(\zeta')
\end{equation} where \begin{equation}\psi(\zeta') = \int_{\C^2} e^{- |a'|^2}
\log|a_1|  \log|\langle a',\zeta'\rangle | da' \quad (\zeta'\in
S^3\subset\C^2)\;.\label{psi}\end{equation}  (To be precise, we have a
well-defined continuous map $\zeta:S^{2d_N-1}\times S^{2d_N-1}\to
S^3/S^1=\C\PP^1$ and $\psi(\zeta')=\psi(\zeta(x,y))$.) By the
Cauchy-Schwartz inequality, we have \begin{eqnarray*} |\psi(\zeta')| &\leq &
\left[\int_{\C^2} e^{- |a'|^2} (\log|a_1|)^2
da'\right]^{1/2}\left[\int_{\C^2} e^{- |a'|^2} (\log|\langle
a',\zeta'\rangle |)^2 da'\right]^{1/2}\\ &=& \int_{\C^2} e^{- |a'|^2}
(\log|a_1|)^2 da'\ =\ C\ <\ +\infty\;,\end{eqnarray*} for all $\zeta'\in
S^3$. Since $$d\sigma =
\sigma(S^{2d_N-1})d\mu_N=\frac{2\pi^{d_N}}{(d_N-1)!}d\mu_N\;,$$ we
have \begin{equation} G_N(x,y)=\frac{1}{\pi^{d_N}}\left(\wt G_N(x,y) -C_N
\right) = \frac{1}{\pi^2}\psi(\zeta') + C'_N\;.\end{equation}  Thus
\begin{equation} E\left(|(\wt Z^N_s,\phi)-(\omega_N,\phi)|^2\right) \leq
\frac{C}{\pi^4 N^2} \sup\|\ddbar\phi\|^2\label{varqed}\end{equation}

\end{proof}

\subsection{Almost everywhere convergence}

We can now complete the proof of Theorem~\ref{zerotheorem} on the
convergence of the zero sets for a random sequence  of sections of increasing
degree, viewed as an element of the probability space ${\cal S} =
\prod_{N=1}^{\infty} SH^0(M, L^{N})$ with the measure $\mu =
\prod_{N=1}^{\infty} \mu_N$. Recall that we identify the unit sphere
$SH^0(M, L^{N})\subset H^0(M, L^{N})$ with the $(2d_N-1)$-sphere
$S^{2d_N-1}\subset \C^{d_N}$  (using the Hermitian inner product described
in \S\ref{notation}); the measure $\mu_N$ is Haar probability measure
on $S^{2d_N-1}$. 

An element in ${\cal S}$ will be denoted ${\bf s} = \{s_N\}$. Since $$|(\wt
Z_{s_N}, \phi)|\leq (\wt Z_{s_N},\omega^{m-1})\|\phi\|_{{\cal C}^0}=
c_1(L)^m\|\phi\|_{{\cal C}^0}\;,$$ by considering a countable ${\cal
C}^0$-dense family of test forms, we need only consider one test form
$\phi$.  By Lemma~\ref{tianscor}, it suffices to show that $$(\wt
Z_{s_N}-\omega_N,\phi)\to 0 \quad\quad \mbox{\rm almost surely}\;.$$  
Consider the random variables \begin{equation} Y_N({\bf s}) =(\wt
Z_{s_N}-\omega_N,\phi)^2\geq 0\;.\end{equation} By Lemma~\ref{variance},
$$\int_{\cal S} Y_N({\bf s}) d\mu({\bf s}) = O(\frac{1}{N^2})\;.$$ Therefore
$$\int_{\cal S}\sum_{N=1}^\infty Y_N d\mu = \sum_{N=1}^\infty \int_{\cal S}
Y_N d\mu <+\infty,$$ and hence $Y_N\to 0$ almost surely. \qed

\medskip\noindent{\it Remark:\/} Since Lemma~\ref{variance} gives an
$O\left(\frac{1}{N^2}\right)$ bound on $E(Y_N)$, we have for any
$\epsilon>0$, $\int_{\cal S} N^{1-2\epsilon} Y_Nd\mu =
O\left(\frac{1}{N^{1+2\epsilon}}\right)$.  Thus the above proof actually
shows that $$\left|\left(\wt Z_N,\phi\right)-(\omega, \phi)\right|\leq
O\left(\frac{1}{N^{\half-\epsilon}}\right)\;, \quad \mbox{\rm almost
surely}.$$

\subsection{Zeros of random orthonormal bases}\label{RONB}

We now switch our attention to 
sequences of orthonormal bases and prove Theorem~\ref{RANTH}. We let
${\cal ONB} = \prod_{N=1}^{\infty} {\cal ONB}_N$ denote the space of
sequences $\{(S_1^N,\ldots,S^N_{d_N}): N = 1, 2, \dots \}$, where
$(S_1^N,\ldots,S^N_{d_N})$ is an element of the space ${\cal ONB}_N$  of
orthonormal bases  for $H^0(M, L^{N})$.  Choosing a fixed $${\bf
e} = \{e_j^N: j = 0, \dots, d_N,\ N=1,2,\ldots \}\in{\cal ONB}$$ gives the
identifications  ${\cal ONB}_N \equiv \U(d_N)$ (the unitary group of rank
$d_N$) and    \begin{equation} {\cal ONB} \equiv \prod_{N=1}^{\infty} 
\U(d_N). \end{equation} We give ${\cal ONB}$ the measure
\begin{equation} \nu:= \prod^{\infty}_{N=1} \nu_N \;,\end{equation} where
$\nu_N$ is the unit-mass Haar measure on $\U(d_N)$.

The variance analogous
of Lemma~\ref{variance} carries over to orthonormal bases:

\begin{lem}\label{variance2} For a smooth test form $\phi$, we have
$$E\big( (\wt
Z^N_{S^N_j}-\omega_N,\phi)^2\big)=O(\frac{1}{N^2})$$ \end{lem}

\begin{proof} Let $\pi^N_j:{\cal ONB}_N\to SH^0(M, L^{N})$ denote
the projection to the $j$-th factor.  Since $\pi^N_{j*} \nu_N
= \mu_N$, we see that
$$E_{\U(d_N)}\big( (\wt
Z^N_{S^N_j}-\omega_N,\phi)^2\big)=E_{S^{2d_N-1}}\big( (\wt
Z^N_s-\omega_N,\phi)^2\big)\;,$$ and thus Lemmas \ref{variance} and
\ref{variance2} are equivalent.\end{proof}

The proof of Theorem~\ref{RANTH} follows easily from
Lemma~\ref{variance2} exactly as in the proof of Theorem~\ref{zerotheorem}.
(The equivalence of the second conclusion follows from \cite[\S 1.3]{Z2}.)

\section{Zeros of $\SU(k)$ polynomials}\label{section-SU}

As an example, we apply Lemma~\ref{mean} to the case
$M=\C\PP^m$, $L=\Oo(1)$, where we give $L$ the standard Hermitian metric
$h_\FS$, whose curvature is the Fubini-Study \kahler form
$\omega=\omega_\FS$ on $\C\PP^m$. We also extend Lemma~\ref{mean} to the
case of simultaneous zeros.

\subsection{$\SU(2)$ polynomials}\label{su(2)}

First consider $m=1$.  Elements of
$H^0(M,L^{N})=H^0(\C\PP^1,\Oo(N))$ are homogeneous polynomials in two
variables of degree $N$, or equivalently, polynomials in one variable of
degree $\leq N$. A basis is given by $\sigma_j=z^j,\ j=0,\ldots,N$.  The inner
product in $H^0(M,L^{N})$ is given by
$$\langle\sigma_j,\sigma_k\rangle = \int_\C\frac{z^j \bar
z^k}{(1+|z|^2)^N}\omega =\frac{1}{\pi}\int_\C\frac{z^j \bar
z^k}{(1+|z|^2)^{N+2}}dxdy\;.$$ Writing the integral in polar coordinates, we
see that the $\sigma_j$ are orthogonal, and \begin{equation}\label{norm}
|\sigma_j|^2=2\int_0^\infty \frac{r^{2j+1}}{(1+r^2)^{N+2}} dr =
\frac{1}{(N+1){{\textstyle{N\choose j}}}}\;.\end{equation}

We thus can choose an orthonormal basis
$$S^N_j= (N+1)^\half{\textstyle{N\choose j}}^\half z^j,\quad\quad
j=0,\ldots,N\;.$$ Next, we note that $$\sum_{j=1}^N \|S^N_j\|^2=(1+|z|^2)^{-N}
\sum_{j=1}^N (N+1){\textstyle{N\choose j}} |z^{2j}| \equiv N+1\;,$$ and thus
$\omega_N=\frac{1}{N}\Phi_N^*\omega_\FS=\omega$.  We thus recover the
following result of \cite[Appendix~C]{B.B.L} on `random $\SU(2)$
polynomials': 

\begin{theo} {\rm \cite{B.B.L}} Suppose we have a random polynomial $$P(z) =
c_0+ c_1 z + \cdots + c_N z^N\;,$$ where ${\rm Re}\,c_0,\ {\rm Im}\,c_0,\
\ldots,\ {\rm Re}\,c_N,\ {\rm Im}\,c_N$ are independent Gaussian random
variables with mean 0 and variances $$E\left(({\rm Re}\,c_j)^2\right)=
E\left(({\rm Im}\,c_j)^2\right)= {\textstyle{N \choose j}}\;.$$ Then the
expected distribution of zeros of $P$ is uniform over $\C\PP^1\approx
S^2$.\label{Bogo}\end{theo}

In fact, Theorem~\ref{zerotheorem} tells us that for a random sequence of
such polynomials, the distribution of zeros approaches uniformity.

\subsection{$\SU(m+1)$ polynomials}\label{su(m)}

We now
turn to the case of polynomials in several variables. An
`$\SU(m+1)$ polynomial of degree $N$' is an element of the probability  space
of homogeneous polynomials of degree $N$ on $\C^{m+1}$ with an $\SU
(m+1)$-invariant Gaussian probability measure. Recall that this
space can be identified with $H^0(\C\PP^m,\Oo(N))$. We give
$H^0(\C\PP^m,\Oo(N))$ the standard inner product. 
A basis for $H^0(\C\PP^m,\Oo(N))$ is given by the monomials
$$\sigma_J=z_0^{j_0}\cdots z_m^{j_m},\quad\quad J=(j_0,\ldots,j_m),\
|J|=N\;.$$  One easily sees that the $\sigma_J$ are orthogonal. 
We compute
\begin{equation}\label{simple-int}|\sigma_J|^2 = 
\int_{\C\PP^m}\frac{|\sigma_J(z)|^2}{|z|^{2N}}\omega^m_{FS}=
\int_{S^{2m+1}}|\sigma_J(z)|^2d\mu^{2m+1} =  \frac{m!j_0!\cdots
j_m!}{(N+m)!}\end{equation}
(where $\mu^{2m+1}$ is Haar probability measure on $S^{2m+1}$),
by writing $$\int_{\C^{m+1}}e^{-|z|^2} |\sigma_J(z)|^2 dz
=\left(\int_\C e^{-|z_0|^2} |z_0|^{2j_0}dz_0\right)\cdots
\left(\int_\C e^{-|z_m|^2} |z_m|^{2j_m}dz_m\right)\;.$$

Therefore, the sections $$S^N_J: =
\left[\frac{(N+m)!}{m!j_0!\cdots j_m!}\right]^\half z^J$$ form an
orthonormal basis for $H^0(\C\PP^m,\Oo(N))$. Furthermore \begin{equation}
\label{sumsq} \sum_{|J|=N}\|S^N_J\|^2
\equiv {\textstyle {N+m\choose m}}\;,\end{equation}
since the sum is $\SU(m+1)$ invariant, hence constant, and the integral of the
left side equals $\dim H^0(\C\PP^m,\Oo(N))$.

In our results on zeros, we can replace the unit sphere $SH^0(M,L^N)$ with
the complex $d_N$-dimensional vector space $H^0(M,L^N)$ with the Gaussian
probability measure $\frac{1}{\pi^{d_N}}e^{-|s|^2} ds$ (where
$ds$ means $2d_N$-dimensional Lebesgue measure).  (We continue to use
the inner product structure on $H^0(M,L^N)$ introduced in
\S\ref{notation}.)  The space of $\SU(m+1)$ polynomials of degree $N$
is by definition the space  $H^0(\C\PP^{m},\Oo(N))$ of homogeneous
polynomials of degree $N$ in $m+1$ variables (or equivalently, polynomials in
$m$ variables of degree \mbox{$\leq N$}) with this Gaussian measure.  We
can use (\ref{simple-int}) to describe the space of $\SU(m+1)$
polynomials explicitly as follows.  For $P\in H^0(\C\PP^m,\Oo(N))$, we write
\begin{equation}\label{gpoly}
P(z_0,\ldots,z_m)=\sum_{|J|=N}\frac{a_J}{\sqrt{j_0!\cdots j_m!}}
z_0^{j_0}\cdots z_m^{j_m}\;.\end{equation}  The Gaussian measure
on $H^0(\C\PP^m,\Oo(N))$ is then given by $$\frac{1}{\pi^{d_N}} e^{-|A|^2}
dA\;, \quad\quad A=(a_J)\in\C^{d_N}\;,$$ where $d_N={N+m \choose m}$.

Lemma 3.1 and (\ref{sumsq}) now tell us that if $P$ is  a polynomial given by
(\ref{gpoly}), with the $a_J$ being independent Gaussian random variables
with mean 0 and variance 1, then the
expected zero current $Z_P$ equals $N\omega_{FS}$.  (This fact, which is the
higher dimensional analogue of Theorem~\ref{Bogo}, is extended to cover
simultaneous zeros in Proposition~\ref{mean3} below.)  Furthermore,
Theorem~\ref{zerotheorem} yields the following:

\begin{prop} Suppose we have a sequence of polynomials
$$P_N(z_0,\ldots,z_m)=\sum_{|J|=N}\frac{a^N_J}{\sqrt{j_0!\cdots j_m!}}
z_0^{j_0}\cdots z_m^{j_m}\;,$$ where the $a^N_J$ are independent Gaussian
random variables with mean 0 and variance 1.  Then 
$$\frac{1}{N}Z_{P_N}\to\omega_{FS}\quad\quad \mbox{\rm almost surely}$$
(weakly in the sense of measures).\end{prop}

\subsection{Expected distribution of simultaneous zeros}

We take a brief detour now to generalize Lemma~\ref{mean} and
Proposition~\ref{mean1} to simultaneous zero sets of holomorphic sections.
This yields a generalization (Proposition~\ref{mean3}) of
Theorem~\ref{Bogo} to the case of simultaneous zeros of polynomials in
several variables. In particular, the 0-dimensional case of
Proposition~\ref{mean3} says that the simultaneous zeros of $m$ random
$\SU(m+1)$ polynomials are uniformly distributed on $\C\PP^m$ with respect
to the volume $\omega^m_\FS$. 

Let
$1\leq \ell\leq m$, and consider the Grassmannian of $\ell$-dimensional
subspaces of $H^0(M,L^{N})$, which we denote $G_\ell
H^0(M,L^{N})$. For an element ${\cal S}={\rm Span}\{s_1,\ldots,s_\ell
\} \in G_\ell H^0(M,L^{N})$, we let $Z_{\cal S}\in{\cal
D}'{}^{\ell,\ell}$ denote the current of integration over the set $\{z\in M:
s_1(z)=\cdots = s_\ell(z)=0\}$. Note that this definition is independent of
the choice of basis $\{s_j\}$ of ${\cal S}$; furthermore by Bertini's
theorem (see \cite{G.H}), the zero sets $Z_{s_j}$ are smooth and intersect
transversely for almost all ${\cal S}$, so we can ignore multiplicities if
we wish.  As before, we consider the normalized current $$\wt Z^N_{\cal S} =
\frac{1}{N^\ell} Z_{\cal S}\;,$$ which we regard as a random current with
${\cal S}$ varying over the probability space $G_\ell H^0(M,L^{N})$
with unit-mass Haar measure. The expected value of $\wt Z^N_{\cal S}$ is
then given by the following elementary formula:

\begin{lem}\label{mean2}  For $N$ sufficiently large so
that $\Phi_N$ is defined, we have:
$$E (\wt Z_{\cal S}^N) =\omega_N^\ell\;.$$
\end{lem}

\begin{proof} Using our fixed orthonormal basis $\{S^N_j\}$, we
can write $s_k=\sum_{j=1}^{d_N}a^j_k S^N_j$.  Let $${\cal S}^\perp =
\{w\in\C\PP^{d_N-1}: \sum_{j=1}^{d_N}a^j_k w_j =0,\ \ k=1,\ldots,\ell\}\;.$$
We let $[{\cal
S}^\perp]$ denote the current of integration over ${\cal S}^\perp$, regarded
as a ${\cal D}'{}^{\ell,\ell}(\C\PP^{d_N-1})$-valued random variable. Since
$\wt Z_{\cal S}^N =\frac{1}{N^\ell}\Phi_N^*  [{\cal S}^\perp]$, we then have
$$E(\wt Z_{\cal S}^N) =\frac{1}{N^\ell}\Phi_N^*  E([{\cal S}^\perp])\;,$$
where $$E([{\cal S}^\perp]) = \int_{G_\ell \C^{d_N}}[{\cal S}^\perp]
d{\cal S}\;.$$ We note that $E([{\cal S}^\perp])$ is $\U(d_N)$-invariant. 
It is well-known that the only $(\ell,\ell)$-currents on projective space
that are invariant under the unitary group are multiples of
$\omega_\FS^\ell$; see \cite[Lemma~3.3]{Sh}. Since $(E([{\cal
S}^\perp]),\omega^{m-\ell})=1$, we conclude that $E([{\cal S}^\perp]) =
\omega^\ell_\FS$ and thus $$E (\wt Z_{\cal S}^N) = \frac{1}{N^\ell}\Phi_N^*
\omega^\ell_\FS =\omega_N^\ell\;.$$ \end{proof}

Applying Corollary~\ref{tianscor}, we obtain the following generalization of
Proposition~\ref{mean1}:

\begin{prop} Let ${\cal S}$ be a random element of $G_\ell H^0(M,L^{N})$,
where $1\leq\ell\leq m$.  Then $$E(\wt Z^N_{\cal S}) = \omega^\ell
+O(\frac{1}{N})\;.$$\end{prop}

We now apply Lemma~\ref{mean2} to random $\SU(m+1)$ polynomials to obtain:

\begin{prop}\label{mean3} Choose  an
$\ell$-tuple ${\cal P} =(P_1,\dots,P_\ell)$ of $\SU(m+1)$ polynomials of
degree $N$ at random. Then
$$E(Z_{\cal P})=N^\ell\omega^\ell_{FS}\;,$$ and in particular
$$E\left({\rm Vol}_{2m-2\ell}\{z\in
U:P_1(z)=\cdots=P_\ell(z)=0\}\right)=\frac{m!}{(m-\ell)!}N^\ell {\rm
Vol}_{2m}(U)$$ for all open subsets $U$ of $\C\PP^m$ (where $\vol_k$ denotes
the Riemannian $k$-volume in $(M,\omega)\,$).\end{prop}

\begin{proof} An  $\ell$-tuple of $\SU(m+1)$  polynomials is an element
of the probability space $$\left(\left[
H^0(\PP^m,\Oo(N)]\right)^\ell, d{\cal G}\right)\;,$$ where $d{\cal G}$ the
$\ell$-fold self-product of the Gaussian measure on
$H^0(\PP^m,\Oo(N))$ (which, of course, is itself a Gaussian measure).  By
(\ref{sumsq}), we conclude as before that $\omega_N=\omega$.  Let $$\Omega =
\left\{(W_1,\ldots,W_\ell)\in \left[ H^0(\PP^m,\Oo(N))\right]^\ell:
W_1\wedge\cdots\wedge W_\ell \neq 0 \right\}\;,$$ and let $\gamma:\Omega
\to G_\ell H^0(\PP^m,\Oo(N))$ be the natural map.  The conclusion follows
from Lemma~\ref{mean2} by noting that $\gamma_*(d{\cal G})$ equals Haar
measure on $G_\ell H^0(\PP^m,\Oo(N))$. \end{proof}

\section{Ergodic orthonormal bases and sections}\label{section-ergodic}

We now turn to the distribution of zeros of sections which form an `ergodic
orthonormal basis'.  As will be explained below, eigenfunctions of quantum
ergodic maps form such a basis.  So do random orthonormal bases.  Both of
these facts belong to now familiar genres of results in quantum chaos.  Let
us briefly recall the basic definitions and results and then prove the
principal new results, Theorem~\ref{QE} and Lemma~\ref{W*}. Proofs of the
background results on ergodic bases are given in the Appendix.

\subsection{The ergodic property}

The weak*-convergence hypothesis of Lemma~\ref{W*}
is closely related to the following `ergodic property':
\medskip

\noindent {\it Definition:\/} We say that ${\bf S} \in {\cal ONB}$  has
the ergodic property  if
$$\lim_{N \rightarrow \infty}\frac{1}{N}
\sum_{n=1}^N \frac{1}{d_n} \sum_{j = 1}^{d_n} \left|\int_M \phi(z) 
\|S^n_j(z)\|^2_{h_n}dV -\bar\phi\right|^2=0\;,\quad \forall \phi \in
{\cal C}(M)\;. \eqno({\cal EP})$$ 
\medskip Here,  $\bar\phi =\int_M \phi dV$ denotes the average value of a
function $f$ over $M$. 

As is well-known (see, for example \cite[\S 1]{Z2}), this property may be
rephrased in the following way: Let ${\bf S} = \{(S^N_1, \dots, S^N_{d_N}):
N=1,2,\ldots\}\in {\cal ONB}$.  Then the ergodic property $({\cal EP})$ is
equivalent to the following weak* convergence property: {\it There exists a
subsequence $\{S'_1,S'_2,\ldots\}$ of relative density one of the sequence
$\{S^1_1, \dots, S^1_{d_1},\ \dots,\ S^N_1,\dots,S^N_{d_N},\ \dots\}$ such
that}
$$\int_M \phi(z) \|S'_n (z) \|^2
dV  \rightarrow \bar\phi\;,\quad \forall \phi \in
{\cal C}(M)\;.
\eqno({\cal EP}')$$ 
A subsequence $\{a_{k_n}\}$ of a sequence
$\{a_n\}$ is said to have relative density one if $\lim_{n\to\infty} n/k_n =
1$.  The equivalence of $({\cal EP})$ and $({\cal EP}')$ is a
consequence of the fact that if
$$\{a_1,a_2,a_3,\dots\}=\{A^1_1,\dots,A^1_{d_1},\
\dots,\ A^n_1,\dots,A^n_{d_n},\ \dots\}$$ is a sequence of non-negative real
numbers, then the following are equivalent: \begin{itemize} \item[\rm i)]
there exists a subsequence $\{a_{k_n}\}$ of relative density one such that
$\lim_{n\to\infty} a_{k_n} \to 0$. \item[\rm ii)] 
$\lim_{N\to\infty}\frac{1}{N}\sum_{n=1}^N a_n \to 0$. \item[\rm iii)]
$\lim_{N\to\infty}\frac{1}{N}\sum_{n=1}^N \frac{1}{d_n}\sum_{j=1}^n A^n_j \to
0$. \end{itemize} 
The equivalence of (i) and (ii) is given in \cite[Theorem~1.20]{Wa}.  (By a
diagonalization argument, one can pick a subsequence independent
of $\phi$ satisfying $({\cal EP}')$.)  For the
equivalence of (ii) and (iii), which depends on the fact that $d_n\sim
n^m$, see \cite[\S 1.3]{Z2}.

We first have:

\begin{theo}\label{EONB} {\bf (a)} A random ${\bf S}\in{\cal ONB}$ has the
ergodic property $({\cal EP})$, or equivalently, $({\cal EP}')$.  In fact,
in complex dimensions $m \geq 2$, a random ${\bf S}\in{\cal ONB}$ has the
property $$\lim_{N \rightarrow \infty} \frac{1}{d_N} \sum_{j=1}^{d_N} \left|
\int_M \phi \|S^N_j\|^2 dV - \bar{\phi}\right|^2 = 0\;,\quad \forall \phi \in
{\cal C}(M)\;,$$ or equivalently, for each $N$ there exists a subset
$\Lambda_N \subset \{1, \dots, d_N\}$ such that $ \frac{\# \Lambda_N}{d_N}
\rightarrow 1$ and $$\lim_{N \rightarrow \infty, j \in \Lambda_N} \int_M \phi
\|S^N_j\|^2 dV = \bar{\phi}.$$\\ {\bf (b)} A random sequence of
sections ${\bf s}=\{s_1,s_2,\dots\}\in{\cal S}$ has a subsequence
$\{s_{N_k}\}$ of relative density 1 such that $$\int_M \phi(z) \|s_{N_k} (z)
\|^2 dV \rightarrow \bar\phi\;,\quad \forall \phi \in {\cal C}(M)\;.$$ In
complex dimensions $m \geq 2$, the entire sequence has this property.
\end{theo} 

Theorem~\ref{EONB}(a) is the line-bundle analogue of
Theorem~(b) in \cite{Z2} on random orthonormal combinations of eigenfunctions
of positive elliptic operators with periodic bicharacteristic flow.  The proof
of Theorem~\ref{EONB} closely parallels those of \cite{Z1,Z2} and strengthens
them in dimensions $m \geq 2$.  Details will be given the Appendix below.

The second setting in which
ergodic orthonormal bases appear is that of quantum ergodicity.  We recall
the following result from \cite[Theorem~B-Corollary~B]{Z3}, which
together with Lemma~\ref{W*} yields Theorem~\ref{QE}. 
 
\begin{theo}\label{Shn} {\rm \cite{Z3}} Let $\{S^N_j\}$ be an
orthonormal basis of eigenfunctions of an ergodic quantum map $U_{\chi, N}$
on $H^0(M, L^{N})$ (as described in Theorem~\ref{QE}).  Then
$\{S^N_j\}$ has the ergodic property $({\cal EP})$, or equivalently, $({\cal
EP}')$. \end{theo}

Theorem~\ref{Shn}
belongs to a long line of results originating in the work of A.~Shnirelman
\cite{Shn.1} in 1974  (see also \cite{Shn.2}) on eigenfunctions of the
Laplacian on compact Riemannian manifolds with ergodic geodesic flow.  The
definition of `quantum map' and the proof of  ergodicity of eigenfunctions
for ergodic  quantum maps over compact \kahler manifolds is contained in
\cite{Z3}, where further references can be found to the literature of
quantum ergodicity. 

We now complete the proofs of Theorems \ref{RANTH} and \ref{QE} by verifying
Lemma~\ref{W*}.

\subsection{Proof of Lemma \ref{W*}}\label{W*proof}

Let $(L,h) \rightarrow (M, \omega)$ and $s_N \in H^0(M, L^{N})$,
$N=1,2,\ldots$, be as in the hypotheses of Lemma~\ref{W*}.  We write
$$u_N =\frac{1}{N} \log \|s_{N}(z)\|_{h_{N}}\;.$$
First we observe that it suffices to show that $u_N\to 0$ in $L^1(M)$. 
Indeed, if that is the case, then for any smooth test form $\phi\in{\cal
D}^{m-1,m-1}(M)$, we have by the Poincar\'e-Lelong formula (\ref{PL}),
$$\left(\frac{1}{N} {Z_N}-\omega,\phi\right)=\left(u_N,
\frac{\sqrt{-1}}{\pi}\ddbar \phi\right)\to 0\;.$$
Since by(\ref{constantmass}), $$\left(\frac{1}{N}{Z_N},\phi\right)\leq
\frac{c_1(L)^m}{m!}\sup|\phi|\;,$$ the conclusion of the lemma holds for
all ${\cal C}^0$ test forms $\phi$.

Next, we observe that: 
\begin{itemize} \item[\rm i)] the functions $u_N$ are
uniformly bounded above on $M$;
\item[\rm ii)] \ $\limsup_{N \rightarrow \infty} u_N\leq 0$.
\end{itemize} Indeed, since $\|s_N\|^2$ converges weakly to 1, we have
$$|s_N|^2=\int_M\|s_N\|^2_{h_N}dV\to 1;.$$ Choose orthonormal
bases $\{S^N_j\}$ and write $s_N=\sum_j a_j S^N_j$, so that
$\sum |a_j|^2 = |s_N|^2$. By Theorem~\ref{tianstheorem}, we have $$\| s_N
(z)\|_{h_N}^2 \leq  |s_N|^2 \sum_{j=1}^{d_N} 
\|S^N_j(z)\|_{h_N}^2 = \left(\frac{c_1(L)^m}{m!} + O(1/N)\right)N^m.$$ Hence
$\| s_N (z)\|_{h_N} \leq C N^{m/2}$ for some $C <\infty$ and taking the
logarithm gives both statements.

Let $e_L$ be a local holomorphic frame for $L$ over $U \subset M$ and let
$e_L^{N}$ be the corresponding frame for $L^{N}.$ Let $g(z) = \|e_L(z)\|_{h}$
so that $\|e_L^{N}(z)\|_{h_N} = g(z)^N$. Then we may write $s_N = f_N e_L^{N}$
with $f_N \in {\cal O}(U)$ and $\|s_N\|_{h_N} = |f_N| g^N$.  It is useful to
consider the function $$v_N = \frac{1}{N}\log |f_N| = u_N -\log g\;,$$ which
is plurisubharmonic on $U$.  (For the properties of plurisubharmonic functions
used here, see for example, \cite{Kl}.)

To finish the proof, we follow the
potential-theoretic approach used by Fornaess and Sibony \cite{F.S} in
their proof of the Brolin-Lyubich theorem on the dynamics of rational
functions.  
Let $U'$ be a relatively compact, open subset of $U$. We must show that
$u_N\to 0$ (or equivalently, $v_N\to -\log g$) in $L^1(U')$.
Suppose on the contrary that $u_N\not\to 0$ in $L^1(U')$.  Then we can
find a subsequence $\{u_{N_k}\}$ with $\|u_{N_k}\|_{L^1(U')} \geq
\delta>0$. By a standard result on subharmonic functions (see
\cite[Theorem~4.1.9]{Ho}), we know that the sequence $\{v_{N_k}\}$ either 
converges uniformly to $-\infty$ on $U'$ or else has a subsequence which is
convergent in $L^1(U')$. Let us now rule out the first possibility. If it
occurred,  there would exist $K > 0$ such that for $k \geq K$,
\begin{equation}\frac{1}{N_k} \log \| s_{N_k}(z)\|_{h_{N_k}} \leq -
1.\label{firstposs}\end{equation} However, (\ref{firstposs}) means that $$
\| s_{N_k}(z)\|_{h_{N_k}}^2 \leq e^{- 2 N_k}\;\;\;\;\forall z \in U'\;, $$
which is inconsistent with the hypothesis that $\| s_{N_k}(z)\|_{h_{N_k}}^2
\rightarrow 1$ in the weak* sense.

Therefore there must exist a subsequence, which we continue to denote  by
$\{v_{N_k}\}$, 
which converges in $L^1(U')$ to some $v
\in L^1(U').$  By
passing if necessary to a further subsequence, we may assume that
$\{v_{N_k} \}$ converges pointwise
almost everywhere  in $U'$ to $v$, and hence
$$v(z) = \limsup_{k \rightarrow \infty}u_{N_k}(z)-\log g\leq
-\log g\;\;\;\;\;\;{\rm (a.e)}\;.$$ Now let 
$$v^*(z):= \limsup_{w \rightarrow z} v(w) \leq -\log g$$ be the
upper-semicontinuous regularization of $v$. Then $v^*$ is plurisubharmonic on
$U'$ and $v^* = v$ almost everywhere.

Since $\|v_{N_k}+\log g\|_{L^1(U')}=\|u_{N_k}\|_{L^1(U')} \geq
\delta>0$, we know that $v^*\not
\equiv -\log g$. Hence, for some $\epsilon
> 0$, the open set $ U_\epsilon = \{z\in U': v^* <-\log g - \epsilon\} $ is
non-empty. Let $U''$ be a non-empty, relatively compact, open subset of $
U_\epsilon$; by Hartogs' Lemma, there exists
a positive integer $K$ such that $v^* \leq -\log g -\epsilon/2$ for $z \in
U'',\ k\geq K$; i.e., \begin{equation} \|s_{N_k}(z)\|_{h_{N_k}}^2 \leq e^{-
\epsilon N_k },\;\;\;\;\;z\in U'', \;\;k \geq K, \end{equation}
which contradicts the weak convergence to $1$. 
\qed

\section{Appendix}

In this Appendix, we give a proof of Theorem~\ref{EONB}, closely following
the proof of Proposition 2.1.4(b) in \cite{Z2}.

To simplify things, we write
\begin{equation}\label{A} A^\phi_{nj}({\bf S}) =\left|\int_M \phi(z) 
\|S^n_j(z)\|^2_{h_n}dV-\bar\phi\right|^2\;.\end{equation} 
In view of the isomorphism $ H^0(M, L^{N})
\cong H^2_N(X)$, we may identify ${\cal ONB}$ with the space of orthonormal
bases of eigenfunctions for the operator $\frac{1}{i} \frac{\partial}
{\partial \theta}$ generating the $S^1$ action on $X$.  Assume without loss
of generality that $\phi$ is real-valued, and consider the Toeplitz operators
$T^\phi_N = \Pi_N M_\phi \Pi_N = \Pi_N M_\phi :H^2_N(X)\to H^2_N(X)$, where
$M_\phi$ is multiplication by the lift of $\phi$ to $X$, and
$\Pi_N:{\cal L}^2(X)\to  H^2_N(X)$ is the orthogonal projection.  Then
$T^\phi_N$ is a self-adjoint operator on $H^2_N(X)$, which can be identified
with a Hermitian $d_N\times d_N$ matrix via the fixed basis {\bf e}.  We
then have \begin{equation} A^\phi_{nj}({\bf S}) = \left|(\phi S^n_j,
S^n_j)-\bar\phi  \right|^2 = \left|(T^\phi_n S^n_j, S^n_j)-\bar\phi 
\right|^2 = \left|(U^*_n T^\phi_n U_n e^n_j, e^n_j)-\bar\phi 
\right|^2,\end{equation}
where ${\bf S} = \{U_N\},\ U_N\in \U(d_N)\equiv{\cal ONB}_N$.

We have
\begin{equation} \bar\phi =  \frac{1}{d_n}\int _M
\sum_{j=1}^{d_n}\|e^n_j\|^2 \phi dV +\int_M
\big[1-\frac{1}{d_n}\sum_{j=1}^{d_n}\|e^n_j\|^2\big]\phi dV =
\frac{1}{d_n}{\rm Tr}\;T^\phi_n +O(\frac{1}{n})\;,\end{equation}   where the
last equality is by Theorem~\ref{tianstheorem}.  Therefore,
\begin{equation}\label{Anj'} A^\phi_{nj}({\bf S})=\wt A^\phi_{nj}({\bf
S})+O(\frac{1}{n})\;,\end{equation} where \begin{equation}\label{Anj}\wt
A^\phi_{nj}({\bf S})=\left|(U^*_n T^\phi_n U_n e^n_j,
e^n_j)-\frac{1}{d_n}{\rm Tr}\;T^\phi_n \right|^2 \;.\end{equation}  (The
bound for the $O(\frac{1}{n})$ term in (\ref{Anj'}) is independent of $\bf
S$.) 

Note that
$iT^\phi_N$ can be identified with an element of the Lie
algebra ${\go u}(d_N)$ of
$\U(d_N)$. Let ${\go t}(d)$ denote the Cartan subalgebra of diagonal
elements in ${\go u}(d)$, and let $\|\cdot\|^2$ denote
the Euclidean inner product on ${\go t}(d)$.  Also let
$$J_d:i{\go u}(d)\rightarrow i{\go t}(d)$$
denote the orthogonal projection (extracting the diagonal).
Finally, let $$\bar J_d (H)=\left(\frac{1}{d}{\rm Tr}\;H\right){\rm
Id}_{d}\;,$$ for Hermitian matrices $H\in i{\go u}(d)$.
(Thus, $H=H^0 +\bar J_d(H)$, with $H^0$ traceless, gives us the
decomposition  ${\go u} (d) =  {\go su} (d)\oplus \R$.)

We introduce the random variables:
$$\begin{array}{l}
Y^\phi_n: {\cal ONB} \to [0,+\infty)\\[6pt]
Y^\phi_n({\bf S}) := \|J_{d_n} (U_n^*T^\phi_n U_n) -
\bar J_{d_n}(T^\phi_n)\|^2\end{array}$$
By (\ref{Anj'}) \begin{equation}\label{Y}\frac{1}{d_n}
Y^\phi_n({\bf S}) =\frac{1}{d_n}\sum_{j=1}^{d_n}\wt A^\phi_{nj}({\bf S})
=\frac{1}{d_n}\sum_{j=1}^{d_n}A^\phi_{nj}({\bf S})
+O(\frac{1}{n})\end{equation} (where the $O(\frac{1}{n})$ term is independent
of {\bf S}). Thus,  $({\cal EP})$ is equivalent to: 
\begin{equation}\label{EP*}
\lim_{N\rightarrow\infty}\frac{1}{N}\sum_{n=1}^N
\frac{1}{d_n} Y^\phi_n ({\bf S}) = 0\;,\quad \forall \phi \in
{\cal C}(M)\;. \end{equation}
The main part of the proof of (\ref{EP*}) is to show the following
asymptotic formula for the expected values of the $Y^\phi_n$.

\begin{lem} \quad $\displaystyle E(Y^\phi_n) =\overline{\phi^2} -
(\bar\phi)^2 +o(1)\;.\label{BGcor}$\end{lem}

\medskip Assume Lemma~\ref{BGcor} for the moment.  The lemma implies that
(\ref{EP*}) holds on the average; i.e.,
\begin{equation}\label{EP*ave}\lim_{N\rightarrow \infty} \frac{1}{N}
\sum_{n=1}^N E\left(\frac{1}{d_n}Y^\phi_n\right)=0\;.\end{equation}
Next we note that
$$ \mbox{\rm Var}\left(\frac{1}{d_n}Y^\phi_n\right)
\leq  \sup \left(\frac{1}{d_n}Y^\phi_n\right)^2
\leq \max_j\; \sup (\wt A^\phi_{nj})^2\;.$$
By (\ref{Anj}),$$ \wt A^\phi_{nj}({\bf S}) \leq 4(U^*_n
T^\phi_n U_n e^n_j,e^n_j)^2\leq 4\sup \phi^2
\;,$$ and therefore 
\begin{equation} \mbox{\rm
Var}\left(\frac{1}{d_n}Y^\phi_n\right)\leq 16 \sup\phi^4
<+\infty\;.\end{equation} Since the variances of the independent random
variables $\frac{1}{d_n}Y^\phi_n$ are bounded, (\ref{EP*}) follows from
(\ref{EP*ave}) and the Kolmogorov strong law of large numbers, which gives
part (a) for general dimensions.  In dimensions $m \geq 2$, we obtain the
improved conclusion as follows: From the fact that $E(\frac{1}{d_N}
Y^{\phi}_N) = O(\frac{1}{N^m})$ it follows that
$E\left(\sum_{N=1}^{\infty}\frac{1}{d_N}Y^{\phi}_N\right)<+\infty$ and thus
$\frac{1}{d_N} Y^{\phi}_N \rightarrow 0$ almost surely when $m \geq 2$. The
quantity we are interested in is $$X_N^{\phi}: = \frac{1}{d_N}
\sum_{j=1}^{d_N} \left| \int_M \phi \|S^N_j\|^2 dV - \bar{\phi} \right|^2 =
\frac{1}{d_N} \sum_{j=1}^{d_N} A^{\phi}_{Nj}.$$ However, by (\ref{Y}),
$$\sup_{{\cal ONB}} |X^{\phi}_N - \frac{1}{d_N} Y^{\phi}_N| =
O(\frac{1}{N}).$$ Hence also $X_N^{\phi} \rightarrow 0$ almost surely.

To verify part (b), we note that since $E(\wt A^\phi_{nj})=E(\wt
A^\phi_{n1})$, for all $j$, it follows from (\ref{Y}) that $E(\wt
A^\phi_{n1})= E(\frac{1}{d_n}Y^\phi_n)$.  Thus,  \begin{equation}
\lim_{N\to\infty}\frac{1}{N}\sum_{n=1}^N \wt
A^\phi_{n1}=0\;,\label{EP1}\end{equation} or equivalently, \begin{equation}
\lim_{N\to\infty}\frac{1}{N}\sum_{n=1}^N
A^\phi_{n1}=0\;.\label{EP*1}\end{equation} Part (b) then follows from
(\ref{EP*1}) exactly as before.

It remains to prove Lemma~\ref{BGcor}. Denote the eigenvalues of $T^\phi_n$
by $\lambda_1,\dots,\lambda_{d_n}$ and write $${\cal
S}_k(\lambda_1,\dots,\lambda_{d_n})=\sum_{j=1}^{d_n} \lambda_j^k\;.$$ 
Note
that  \begin{equation} {\rm Tr}\;(T^\phi_n)^k ={\cal
S}_k(\lambda_1,\dots,\lambda_{d_n})\;.
\end{equation}  
We shall use the following
`Szego limit theorem' due to Boutet de Monvel and Guillemin
\cite[Theorem~13.13]{B.G}:

\begin{lem} {\rm \cite{B.G}} For $k\in\Z^+$, we have
$$ \lim_{N \rightarrow \infty}
\frac{1}{d_N} {\rm Tr}\;(T^\phi_N)^k = \overline{\phi^k}\;.$$
\label{BGlemma}\end{lem} 
Lemma~\ref{BGcor} is an immediate consequence of
Lemma~\ref{BGlemma} and the following formula:

\begin{equation}\label{orbit-int} \int_{\U(d)} \|J_d (U^*D(\vec\lambda) U) -
\bar J_d(D(\vec\lambda))\|^2 dU =
\frac{{\cal S}_2(\vec\lambda)}{d+1}-
\frac{{\cal S}_1(\vec\lambda)^2}{d(d+1)}\;,\end{equation}
where $\vec\lambda=(\lambda_1,\dots,\lambda_d)\in\R^d$,
$D(\vec\lambda)$ denotes the diagonal matrix with entries equal to the
$\lambda_j$, and integration is
with respect to Haar probability measure on $\U(d)$.

A proof of the identity (\ref{orbit-int}) is given in \cite[pp.~68--69]{Z1}
(see also \cite{Z2}).  For completeness, we provide here a simplified proof
of (\ref{orbit-int}) following the methods of \cite{Z1,Z2}.  Let ${\cal
E}(\vec\lambda)$ denote the left side of (\ref{orbit-int}). Since ${\cal
E}(\vec\lambda)$ is a homogeneous, degree 2, symmetric polynomial in
$\vec\lambda$, we can write \begin{equation}\label{solve}{\cal
E}(\vec\lambda)=c_d {\cal S}_2(\vec\lambda)+c'_d {\cal S}_1(\vec\lambda)^2
\;.\end{equation}  Substituting $\vec\lambda = (1,\dots,1)$ in
(\ref{solve}) and using the fact that ${\cal E}(1,\dots,1)=0$, we  conclude
that $c'_d= -c_d/d$.  To find $c_d$, we substitute  $\vec\lambda =
(1,0,\dots,0)$, and write $D=D(1,0,\dots,0)$.  For $U=(u_{jk})\in \U(d)$,
we have $$(U^*DU)_{jj}=|u_{1j}|^2\;,\quad \bar J_dD=\frac{1}{d}{\rm
Id}_d\;.$$ Therefore, \begin{eqnarray*}{\cal E}(1,0,\dots,0) &=&
\int_{\U(d)}\sum_{j=1}^d\left(|u_{1j}|^2-\frac{1}{d}\right)^2 dU =
\int_{S^{2d-1}}\sum_{j=1}^d\left(|a_j|^2-\frac{1}{d}\right)^2
d\mu^{2d-1}(a)\\
&=&\int_{S^{2d-1}}\big(\sum_{j=1}^d|a_j|^4-\frac{1}{d}\big) d\mu^{2d-1}(a)
= -\frac{1}{d} + d\int_{S^{2d-1}}|a_1|^4d\mu^{2d-1}(a)\;, \end{eqnarray*}
where $a=(a_1,\dots,a_d)\in S^{2d-1}$ and $\mu^{2d-1}$ is unit-mass Haar
measure on $S^{2d-1}$. By (\ref{simple-int}),
$$\int_{S^{2d-1}}|a_1|^4d\mu^{2d-1}(a)=\frac{2}{d(d+1)}\;,$$ and therefore
\begin{equation}\label{E} {\cal E}(1,0,\dots,0) =
\frac{d-1}{d(d+1)}\;.\end{equation}  Substituting (\ref{E}) into
(\ref{solve}) with $c'_d= -c_d/d$, we conclude that
$$c_d=\frac{1}{d+1}\;.$$ \qed

\end{document}